\documentclass[12pt,leqno]{amsart}
\newtheorem{thm}{Theorem}[section]
\newtheorem{defi}{Definition}[section]
\newtheorem{lem}{Lemma}[section]
\newtheorem{cor}{Corollary}[section]

\newtheorem{rem}{Remark}[section]

\numberwithin[section]{equation}{section}
\newcommand{\be}{\begin{equation}}
\newcommand{\ee}{\end{equation}}
\numberwithin{equation}{section}
\newcommand{\bea}{\begin{eqnarray}}
\newcommand{\eea}{\end{eqnarray}}
\newcommand{\beb}{\begin{eqnarray*}}
\newcommand{\eeb}{\end{eqnarray*}}
\usepackage{amssymb,amsfonts,amsthm,setspace,indentfirst}
\numberwithin{equation}{section}
\usepackage[dvips]{graphics}
\usepackage{epsfig}
\begin{document}
\title[On pseudosymmetric manifolds]{\bf{On pseudosymmetric manifolds}}
\author[A. A. Shaikh, R. Deszcz, M. Hotlo\'{s}, J. Je\l owicki and H. Kundu]{Absos Ali Shaikh$^*$, 
Ryszard Deszcz, Marian Hotlo\'{s}, Jan Je\l owicki and Haradhan Kundu}
\date{}
\address{\newline\noindent Absos Ali Shaikh and Haradhan Kundu,\newline Department of Mathematics,\newline University of 
Burdwan, Golapbag,\newline Burdwan-713104,\newline West Bengal, India}
\email{aask2003@yahoo.co.in, \ aashaikh@mathburuniv.ac.in, \ kundu.haradhan@gmail.com}
\address{\noindent\newline Ryszard Deszcz and Jan Je\l owicki, \newline Department of Mathematics,\newline Wroc\l aw University of Environmental and Life Sciences\newline Grunwaldzka 53, 50-357 Wroc\l aw , Poland}
\email{ryszard.deszcz@up.wroc.pl \ jan.jelowicki@up.wroc.pl}
\address{\noindent\newline Marian Hotlo\'{s}\newline Institute of Mathematics and Computer Science\newline Wroc\l aw 
University of Technology\newline  Wybrze\.{z}e Wyspia\'{n}skiego 27, 50-370 Wroc\l aw, Poland}
\email{marian.hotlos@pwr.edu.pl}
%
\begin{abstract}
In the literature, there are two different notions of pseudosymmetric manifolds, 
one by Chaki \cite{CHA} and other by Deszcz \cite{DES}, and there are many papers related to these notions. 
The object of the present paper is to deduce necessary 
and sufficient conditions for a Chaki pseudosymmetric \cite{CHA} (resp. pseudo Ricci symmetric \cite{CHA1}) manifold 
to be Deszcz pseudosymmetric (resp. Ricci pseudosymmetric). We also study the necessary and sufficient conditions 
for a weakly symmetric \cite{tb} (resp. weakly Ricci symmetric \cite{ts}) manifold by Tam$\acute{\mbox{a}}$ssy and Binh 
to be Deszcz pseudosymmetric (resp. Ricci pseudosymmetric). 
We also obtain the reduced form of the defining condition of weakly Ricci symmetric manifolds 
by Tam$\acute{\mbox{a}}$ssy and Binh \cite{ts}. Finally we give some examples 
to show the independent existence of such types of pseudosymmetry which also ensure the existence of Roter type 
and generalized Roter type manifolds and the manifolds with recurrent curvature $2$-form (\cite{besse}, \cite{LR}) 
associated to various curvature tensors.
\end{abstract}
\noindent\footnotetext{ $^*$ Corresponding author.\\
$\mathbf{2010}$\hspace{5pt}Mathematics\; Subject\; Classification: 53B20, 53B30, 53C15, 53C21, 53C25. }
\maketitle
\begin{center}
{\sl Dedicated to the Memory of Professor M.C. Chaki}
\end{center}
\section{\bf{Introduction}}\label{intro}
The geometry of a space mainly depends on the curvature of the space. One of the most 
important geometric property of a space is symmetry. The study of symmetry of a manifold 
began with the works of Cartan \cite{ca} and then his notion has been weakened by various authors in different 
directions with several defining conditions by giving some curvature restrictions. Cartan first 
classified complete simply connected locally symmetric spaces \cite{ca} for the Riemannian case and the same was done for non-Riemannian case by Cahen and Parker (\cite{CAH}, \cite{CAH1}). Later various weaker symmetries are studied  as generalizations or extensions  of Cartan's notion such as recurrent manifolds by Walker \cite{Ag}, generalized recurrent manifolds by Dubey \cite{DUB}, quasi-generalized recurrent manifolds by Shaikh and Roy \cite{ROY1}, weakly generalized recurrent manifolds by 
Shaikh and Roy \cite{ROY}, hyper-generalized recurrent manifolds by Shaikh and Patra \cite{SP}, semisymmetric manifolds by Cartan \cite{ca} 
(and classified in the Riemannian case by Szab$\acute{\mbox{o}}$ \cite{sz}), 
pseudosymmetric manifolds by Deszcz \cite{DES}, pseudosymmetric manifolds by Chaki \cite{CHA}, weakly symmetric manifolds by Selberg \cite{sel}, weakly symmetric manifolds by Tam$\acute{\mbox{a}}$ssy and Binh \cite{tb} etc. We note that the notion of pseudosymmetry by Deszcz \cite{DES}  is different to that by Chaki \cite{CHA}. Also, the notion of weakly symmetric manifold by Selberg \cite{sel} is different to that by Tam$\acute{\mbox{a}}$ssy and Binh \cite{tb} and throughout the paper we will confine ourselves with the notion of weakly symmetric manifold by Tam$\acute{\mbox{a}}$ssy and Binh \cite{tb}. On the analogy, various types of symmetry, recurrency, weak symmetry and pseudosymmetry are also studied by many authors for other generalized curvature tensors.\\
\indent The notion of Ricci symmetry has been weakened by many authors in several directions such as Ricci recurrent manifolds by Patterson \cite{pat}, Ricci semisymmetric manifolds also named as Ric-semisymmetric manifolds \cite{Lumiste}, Ricci pseudosymmetric manifolds by Deszcz (\cite{adm}, \cite{DES1}), pseudo Ricci symmetric manifolds by Chaki \cite{CHA1}, weakly Ricci symmetric manifolds 
by Tam$\acute{\mbox{a}}$ssy and Binh \cite{ts} etc. 
We note that pseudo Ricci symmetry by Chaki and Ricci pseudosymmetry by Deszcz are also different.\\
\indent Ewert-Krzemieniewski (\cite{sti}, \cite{sti1}) simplified the defining condition of a weakly symmetric manifold 
by Tam$\acute{\mbox{a}}$ssy and Binh and showed that such a manifold and Chaki pseudosymmetric manifold are the same. 
Motivated by the study of Ewert-Krzemieniewski (\cite{sti}, \cite{sti1}), we present three types of weak symmetry, viz., 
type I, II, III and prove that all the three types of weak symmetry are equivalent for a generalized curvature tensor. 
Moreover, for a proper generalized curvature tensor all the three types of weak symmetry are equivalent with Chaki pseudosymmetry. Recently, Mantica and Molinari \cite{mmwsn} studied weakly $Z$-symmetric manifold by assuming $Z$ as a symmetric $(0,2)$-tensor with a specific form and obtained the reduced form of the defining condition of such notion. In \cite{ManSuh12a} the authors defined and investigated pseudo $Z$-symmetric Riemannian manifolds where the tensor $Z$ is of the above form. In the present paper considering $Z$ as a $(0,2)$-tensor (symmetric or skew-symmetric) we obtain the reduced form of weakly $Z$-symmetric manifold (see Theorem \ref{thm3.1}), which entails the simplified form of defining condition of weakly Ricci symmetric manifold. 
The paper is organized as follows. Section 2 deals with preliminaries where the defining conditions of pseudosymmetry 
and weak symmetry are presented. Section 3 is concerned with pseudosymmetry by Chaki and three types of weak symmetry 
for a $(0,k)$-tensor field along with their interrelationship. Although in several papers it was mentioned 
that both the pseudosymmetry by Chaki and Deszcz are different but there are only a few papers where studied among their interrelationship (see \cite{DES}, \cite{MD14}, \cite{ManMol14}, \cite{ManSuh12a}, \cite{ManSuh13}). 
In Section 4 we establish the necessary and sufficient conditions for a weakly symmetric (resp. weakly Ricci symmetric) manifold by Tam$\acute{\mbox{a}}$ssy and Binh to be Deszcz pseudosymmetric (resp. Ricci pseudosymmetric), from which we also deduce the necessary and sufficient conditions for a Chaki pseudosymmetric (resp. pseudo Ricci symmetric) 
manifold to be Deszcz pseudosymmetric (resp. Ricci pseudosymmetric). Finally, in the last section, the independent existence of Chaki pseudosymmetric and Deszcz pseudosymmetric manifolds are given by some proper examples with various new metrics, which compelled us to introduce the notion of generalized Roter type manifolds.
\section{\bf{Preliminaries}}\label{pre}
%
Let $(M,g)$, $n = \dim M \geqslant 3$, be a semi-Riemannian manifold,
i.e. connected smooth manifold equipped with a semi-Riemannian metric $g$. 
We denote by $\nabla$, $R$, $S$, $\kappa$, the Levi-Civita connection, the Riemann-Christoffel curvature tensor, Ricci tensor and scalar curvature of $(M,g)$, respectively. Now for $(0,2)$-tensors $A$ and $D$, their Kulkarni-Nomizu product 
(see, e.g., \cite{GLOG}, \cite{SKe}) $A\wedge D$ is given by
\bea\label{eq2.1}
(A \wedge D)(X_1,X_2,Y_1,Y_2)&=&A(X_1,Y_2)D(X_2,Y_1) + A(X_2,Y_1)D(X_1,Y_2)\\\nonumber
&&-A(X_1,Y_1)D(X_2,Y_2) - A(X_2,Y_2)D(X_1,Y_1),
\eea
where $X_1, X_2, Y_1, Y_2\in \chi(M)$, $\chi(M)$ being the Lie algebra of all smooth vector fields on M. 
Throughout the paper we consider $X, Y, X_i, Y_i \in \chi(M)$, $i = 1,2, ...  $\\
A tensor $B$ of type $(1,3)$ on $M$ is said to be generalized curvature tensor (see, e.g., \cite{SKe}), if
\beb
&(i)&B(X_1,X_2)X_3+B(X_2,X_3)X_1+B(X_3,X_1)X_2=0,\\
&(ii)&B(X_1,X_2)X_3+B(X_2,X_1)X_3=0,\\
&(iii)&B(X_1,X_2,X_3,X_4)=B(X_3,X_4,X_1,X_2),
\eeb
where $B(X_1,X_2,X_3,X_4)=g(B(X_1,X_2)X_3,X_4)$. Throughout the paper we consider $X, Y, X_i, Y_i \in \chi(M)$, $i = 1,2, \cdots $. Here we use the same symbol $B$ for the generalized curvature tensor of type $(1,3)$ and $(0,4)$. Moreover if $B$ satisfies the second Bianchi identity i.e.,
$$(\nabla_{X_1}B)(X_2,X_3)X_4+(\nabla_{X_2}B)(X_3,X_1)X_4+(\nabla_{X_3}B)(X_1,X_2)X_4=0,$$
then $B$ is called a proper generalized curvature tensor. Throughout this paper we denote 
by $B$ the generalized curvature tensor unless otherwise stated.\\
Now for a generalized curvature tensor $B$ and given two vector fields $X,Y\in\chi(M)$ one can define an endomorphism $\mathcal{B}(X,Y)$ by
$$\mathcal{B}(X,Y)(X_1)=B(X,Y)X_1, \ \mbox{for all $X_1\in\chi(M)$}.$$
Again if $X,Y\in\chi(M)$ then for a $(0,2)$-tensor $A$ one can define an endomorphism $X \wedge_A Y$, by
$$(X \wedge_A Y)X_1 = A(Y,X_1)X - A(X,X_1)Y, \ \mbox{for all $X_1 \in \chi(M) $}.$$
Some most useful generalized curvature tensors are Gaussian curvature tensor $G$, Weyl conformal curvature tensor $C$, 
concircular curvature tensor $K$ and conharmonic curvature tensor $conh(R)$, which are respectively given by
$\frac{1}{2} g \wedge g$, 
$R -\frac{1}{n-2} g \wedge S+\frac{\kappa}{2(n-1)(n-2)}g \wedge g$, 
$R - \frac{\kappa}{2 n(n-1)}g \wedge g$, 
$R - \frac{1}{n-2} g \wedge S$. 
We note that the Weyl projective curvature tensor $P$ given by $\mathcal P(X,Y) = \mathcal R(X,Y) -\frac{1}{n-2} X\wedge_S Y$ is not a generalized curvature tensor.\\
Let us now consider $\mathcal T^r_k(M)$ be the space of all tensor fields of type $(r,k)$ on $M$, 
$r,k\in \mathbb N\cup\left\{0\right\}$. Now for $T\in \mathcal T^0_k(M)$, $k\geq 2$, 
and a generalized curvature tensor $B$ one can define a $(0,k+2)$ 
tensor $B\cdot T$ given by (\cite{DDVY}, \cite{SKe})
\beb\label{rdot}
&&B\cdot T(X_1,X_2, \ldots ,X_k;X,Y) = (\mathcal{B}(X,Y) T)(X_1,X_2, \ldots ,X_k)\\\nonumber
&&= -T(\mathcal{B}(X,Y)X_1,X_2, \ldots ,X_k) - \cdots - T(X_1,X_2, \ldots ,\mathcal{B}(X,Y)X_k),
\eeb
and for a $(0,2)$-tensor $A$ one can define a $(0,k+2)$-tensor $Q(A,T)$ 
as (\cite{DDVY}, \cite{SKe}, \cite{tac})
\beb\label{qgr}
&&Q(A,T)(X_1,X_2, \ldots ,X_k;X,Y) = ((X \wedge_A Y) T)(X_1,X_2, \ldots ,X_k)\\\nonumber
&&= -T((X \wedge_A Y)X_1,X_2, \ldots ,X_k) - \cdots - T(X_1,X_2, \ldots ,(X \wedge_A Y)X_k).
\eeb
For an $1$-form $\mu$ and a vector field $X$ on $M$, we can define an endomorphism $\mu_{_X}$ as
\begin{eqnarray*}
\mu_{_X}(X_1) &=& \mu(X_1)X, \ \mbox{for all} \ X_1\in \chi(M).
\end{eqnarray*}
Then we can define $\mu_{_X}$ \cite{SKe} as an operation on a $(0,k)$-tensor field $T$ as follows:
\beb\label{pidot}
&&(\mu_{_X} \cdot T)(X_1,X_2, \ldots ,X_k)\\\nonumber
&&= -T(\mu_{_X}(X_1),X_2, \ldots ,X_k) - \cdots - T(X_1,X_2, \ldots ,\mu_{_X}(X_k)),\\\nonumber
&&= -\mu(X_1)T(X,X_2, \ldots ,X_k) -\mu(X_2)T(X_1,X, \ldots ,X_k)- \cdots -\mu(X_k)T(X_1,X_2, \ldots ,X),
\eeb
for all $X, X_i \in \chi(M)$.\\
\indent For a complete classification with generalized curvature tensor and equivalency of various types of pseudosymmetric conditions, we refer to \cite{SKe} and references therein.
\begin{defi}$($\cite{ca}, \cite{sz}$)$ 
A semi-Riemannian manifold $(M,g)$, $n \geqslant 3$, 
admitting a $(0,k)$-tensor field $T$ 
is said to be $T$-semisymmetric 
if $R\cdot T = 0$ on $M$.
\end{defi}
In particular, if $T= R$ (resp., $S$) then the manifold is called semisymmetric 
(resp., Ricci semisymmetric).
\begin{defi}$($\cite{adm}, \cite{DEF}, \cite{DES}, \cite{des7}$)$ 
A semi-Riemannian manifold $(M,g)$, $n \geqslant 3$, admitting a $(0,k)$-tensor field $T$ 
is said to be Deszcz $T$-pseudosymmetric $($resp., Ricci generalized $T$-pseudosymmetric$)$
if $R \cdot T$ and $Q(g,T)$ $($resp., $R \cdot T$ and $Q(S,T)$ $)$ are linearly dependent, 
i.e., $R \cdot T = L_T\, Q(g,T)$ $($resp., $R \cdot T = L_G\, Q(S,T)$ $)$ holds on the set 
$U_T = \{x\in M : Q(g,T) \neq 0 \ at \ x\}$ $($resp., $U_G = \{x\in M : Q(S,T) \neq 0 \ at \ x\}$ $)$, 
where $L_T$ $($resp., $L_G$ $)$ is some function on $U_T$ $($resp., $U_G$ $)$.
\end{defi}
In particular, if $R \cdot R = L_R\, Q(g,R)$ (resp., $R \cdot S = L_S\, Q(g,S)$) then the manifold is called Deszcz pseudosymmetric (resp., Ricci pseudosymmetric). Throughout the paper we denote Deszcz pseudosymmetric manifold by $(DPS)_n$ and Deszcz Ricci pseudosymmetric manifold by $(DRPS)_n$.
For details about the Deszcz pseudosymmetry, Ricci generalized pseudosymmetry, 
as well other conditions of pseudosymmetry type we refer the reader the papers: 
\cite{adm}, \cite{DEF}-\cite{kun}, \cite{Saw3} and also references therein. We note that \cite{des7} is the first paper, 
in which manifolds satisfying $R \cdot R = L_R\, Q(g,R)$ were called pseudosymmetric manifolds. 
It seems that the Schwarzschild spacetime, the Kottler spacetime, the Reissner-Nordstr\"{o}m spacetime, 
as well as some Friedmann-Lema{\^{\i}}tre-Robertson-Walker spacetimes are the ``oldest'' examples 
of non-semisymmetric pseudosymmetric warped product manifolds (see, e.g., \cite{DESZ}).

\begin{defi}
A semi-Riemannian manifold $(M,g)$, $n \geqslant 3$, 
admitting a $(0,4)$-tensor field $T$ 
is said to be weakly $T$-symmetric by Tam$\acute{\mbox{a}}$ssy and Binh \cite{tb} if
\bea\label{wsntb}
(\nabla_X T)(X_1,X_2,X_3,X_4) &=& \alpha(X)\, T(X_1,X_2,X_3,X_4)\\\nonumber
&+& \beta(X_1)\, T(X,X_2,X_3,X_4) + \bar\beta(X_2)\, T(X_1,X,X_3,X_4)\\\nonumber
&+& \gamma(X_3)\, T(X_1,X_2,X,X_4) + \bar\gamma(X_4)\, T(X_1,X_2,X_3,X)
\eea
holds on the set $U_J = \{x\in M : \nabla T - \xi\otimes T \neq 0 \ \mbox{for any 1-form $\xi$ at} \ x\}$, 
where $\alpha, \beta, \bar\beta, \gamma$ and $\bar\gamma$ are associated $1$-forms and we say that 
($\alpha$, $\beta$, $\bar\beta$, $\gamma$, $\bar\gamma$) is a solution of this weakly $T$-symmetric manifold. 
Especially, if the solution is of the form ($2 \phi$, $\phi$, $\phi$, $\phi$, $\phi$) i.e.,
\begin{eqnarray}\label{cpsn}
(\nabla_{X} T)(X_1,X_2,X_3,X_4)&=& 2\phi(X)\, T(X_1,X_2,X_3,X_4) \\\nonumber
&+& \phi(X_1)\, T(X,X_2,X_3,X_4)+\phi(X_2)\, T(X_1,X,X_3,X_4) \\\nonumber
&+& \phi(X_3)\, T(X_1,X_2,X,X_4)+\phi(X_4)\, T(X_1,X_2,X_3,X)
\end{eqnarray}
holds on the set $U_L = \{x\in M : \nabla T \neq 0 \ at \ x\}$, then the manifold is called Chaki $T$-pseudosymmetric \cite{CHA}. Again if the solution is of the form ($\pi$, $0$, $0$, $0$, $0$) i.e.,
\be\label{kn}
(\nabla_{X} T)(X_1,X_2,X_3,X_4) = \pi(X)\, T(X_1,X_2,X_3,X_4)
\ee
holds on the set $U_L = \{x\in M : \nabla T \neq 0 \ at \ x\}$, then the manifold is called $T$-recurrent \cite{Ag}.
\end{defi}
In particular, if $T= R$ then the manifold satisfying \eqref{wsntb} (resp., \eqref{cpsn}, \eqref{kn}) is called weakly symmetric manifold by Tam$\acute{\mbox{a}}$ssy and Binh (resp., Chaki pseudosymmetric manifold, recurrent). 
Throughout the paper we denote weakly symmetric manifold by $(WS)_n$, Chaki pseudosymmetric manifold by $(CPS)_n$ and recurrent manifold by $K_n$. 
For details about the study of weak symmetry with various curvature tensors and structures, 
we refer the reader the papers \cite{DE}, \cite{DB}, \cite{JANA3} and also references therein. 
For decomposable and warped product weakly symmetric manifolds, we refer the reader to see \cite{Bin} and \cite{SK}.
\begin{defi}
A semi-Riemannian manifold $(M,g)$, $n \geqslant 3$, admitting a $(0,2)$-tensor field $Z$ 
is said to be weakly $Z$-symmetric by Tam$\acute{\mbox{a}}$ssy and Binh \cite{ts} if
\bea\label{wrsntb}
(\nabla_X Z)(X_1,X_2) &=& \delta(X)\, Z(X_1,X_2) + \eta(X_1)\, Z(X,X_2) + \lambda(X_2)\, Z(X_1,X)
\eea
holds on the set $U_Q = \{x\in M : \nabla Z - \xi\otimes Z \neq 0 \ \mbox{for any $1$-form $\xi$ at} \ x\}$, 
where $\delta, \eta$ and $\lambda$ are associated $1$-forms and we say that $(\delta$, $\eta$, $\lambda)$ is a solution 
of this weakly $Z$-symmetric manifold. Especially, if the solution is of the form $(2\psi$, $\psi$, $\psi)$, i.e.,
\be\label{cprsn}
(\nabla_{X} Z)(X_1,X_2)= 2\psi(X)\, Z(X_1,X_2)+\psi(X_1)\, Z(X,X_2)+\psi(X_2)\, Z(X_1,X)
\ee
holds on the set $U_F = \{x\in M : \nabla Z \neq 0 \ at \ x\}$, then the manifold is called Chaki pseudo $Z$-symmetric \cite{CHA1}. Again if the solution is of the form $(\pi$, $0$, $0)$, i.e.,
\be\label{rkn}
(\nabla_{X} Z)(X_1,X_2)= \pi(X)\, Z(X_1,X_2)
\ee
holds on the set $U_F = \{x\in M : \nabla Z \neq 0 \ at \ x\}$, then the manifold is called $Z$-recurrent \cite{pat}.
\end{defi}
A manifold satisfying \eqref{wrsntb} (resp., \eqref{cprsn}) with $Z = S+a g$, $a$ being an arbitrary scalar function, was investigated in \cite{mmwsn} (resp., \cite{ManSuh12a}). In particular, if $Z=S$ then the manifold satisfying \eqref{wrsntb} (resp., \eqref{cprsn}, \eqref{rkn}) is called weakly Ricci symmetric (resp., Chaki pseudo Ricci symmetric, Ricci recurrent) and throughout the paper we denote such a manifold by $(WRS)_n$ (resp., $(CPRS)_n$, $RK_n$). 
For details about the study of $(WRS)_n$, its generalization and related works, we refer the reader the papers 
\cite{DEGS}, \cite{mmid}, \cite{mmwsn}, \cite{S5} and also references therein.  Also for details about the study of $(CPRS)_n$, we refer the reader the papers \cite{CHA1} and also references therein.\\
\indent Recently Mantica and Suh (\cite{ManSuh}, \cite{ManSuh01}, \cite{ManSuh02} and \cite{ManSuh12b}) presented a curvature restriction which is necessary and sufficient 
for the recurrency of a specific curvature $2$-form associated to that curvature tensor. For a generalized curvature tensor $B$, 
the associated $2$-form is defined as (\cite{besse}, \cite{LR})
$$\Omega_{(B)l}^m = B_{jkl}^m dx^j \wedge dx^k,$$
where $\wedge$ indicates the exterior product. Again for a symmetric $(0,2)$-tensor $Z$, the associated $1$-form \cite{SKP} is defined as
$$\Lambda_{(Z)l} = Z_{lm} dx^m.$$
In \cite{ManSuh}, \cite{ManSuh01}, \cite{ManSuh02} and \cite{ManSuh12b} Mantica and Suh showed that $\Omega_{(B)l}^m$ is recurrent (i.e., $\mathcal D \Omega_{(B)l}^m  = \alpha \wedge \Omega_{(B)l}^m$, $\mathcal D$ is the exterior derivative and $\alpha$ is the associated $1$-form) if and only if
$$\nabla _{h} B_{ijkl} +\nabla _{i} B_{jhkl} +\nabla _{j} B_{hikl} = \alpha_{h}\, B_{ijkl} +\alpha_{i}\, B_{jhkl} +\alpha_{j}\, B_{hikl}$$
and $\Lambda_{(Z)l}$ is recurrent (i.e., $\mathcal D \Lambda_{(Z)l}  = \alpha \wedge \Lambda_{(Z)l}$) if and only if
$$\nabla _{i} Z_{kl} - \nabla _{k} Z_{il} = \alpha_{i}\, Z_{kl} - \alpha_{k}\, Z_{il}.$$
In this connection we also note that the curvature restriction
$$\nabla _{i} S_{jk} +\nabla _{j} S_{ki} +\nabla _{k} S_{ij} = \alpha_{i}\, S_{jk} + \alpha_{j}\, S_{ki} + \alpha_{k}\, S_{ij}$$
was investigated by Shaikh and Jana \cite{S5}.\\
Hence for a $(0,4)$-tensor $T$ and a $(0,2)$-tensor $Z$, we can have the following curvature restrictions:
\begin{eqnarray}\label{b1}
\nabla _{h} T_{ijkl} +\nabla _{i} T_{jhkl} +\nabla _{j} T_{hikl} &=& 0
\end{eqnarray}
\begin{eqnarray}\label{b2}
\alpha_{h}\, T_{ijkl} +\alpha_{i}\, T_{jhkl} +\alpha_{j}\, T_{hikl} &=& 0,
\end{eqnarray}
\begin{eqnarray}\label{b3}
\nabla _{h} T_{ijkl} +\nabla _{i} T_{jhkl} +\nabla _{j} T_{hikl} &=&
\alpha_{h}\, T_{ijkl} +\alpha_{i}\, T_{jhkl} +\alpha_{j}\, T_{hikl}.
\end{eqnarray}
\begin{eqnarray}\label{b4}
\nabla _{i} Z_{kl} - \nabla _{k} Z_{il} = \alpha_{i}\, Z_{kl} - \alpha_{k}\, Z_{il},
\end{eqnarray}
where $\alpha$ is the corresponding $1$-form for the restrictions. In Section \ref{Example}, 
we examine the above curvature restrictions for the existence of the manifolds with recurrent curvature $2$-form.
\section{\bf{Weak symmetry and Chaki pseudosymmetry}}\label{wsnpsn}
We consider a weakly $B$-symmetric manifold whose defining condition is given in (\ref{wsntb}) for $T=B$. 
We note that in 1995 Prvanovi\'{c} \cite{pra} showed that in a $(WS)_n$, $\beta = \bar\beta$ and $\gamma = \bar \gamma$ 
and then in 1999 the same was again proved by De and Bandyopadhyay \cite{DB}. 
Hence for a weakly $B$-symmetric manifold the solution ($\alpha$, $\beta$, $\bar\beta$, $\gamma$, $\bar\gamma$) 
turns into ($\alpha$, $\beta$, $\beta$, $\gamma$, $\gamma$). Thus the defining condition of a weakly $B$-symmetric manifold takes the form
\bea\label{eq3.1}
(\nabla_X B)(X_1,X_2,X_3,X_4) &=& \alpha(X)\, B(X_1,X_2,X_3,X_4)\\\nonumber
&+& \beta(X_1)\, B(X,X_2,X_3,X_4) + \beta(X_2)\, B(X_1,X,X_3,X_4)\\\nonumber
&+& \gamma(X_3)\, B(X_1,X_2,X,X_4) + \gamma(X_4)\, B(X_1,X_2,X_3,X).
\eea
Again Ewert-Krzemieniewski (\cite{sti}, \cite{sti1}) proved that in a weakly $B$-symmetric manifold with solution 
($\alpha$, $\beta$, $\bar\beta$, $\gamma$, $\bar\gamma$) there exists another solution ($\alpha$, $\sigma$, $\sigma$, $\sigma$, $\sigma$). 
We note that we can determine the $1$-form $\sigma$ as $\sigma = \frac{\beta+\gamma}{2}$. 
Moreover, if $B$ is a proper generalized curvature tensor then there exists a solution 
($2\epsilon$, $\epsilon$, $\epsilon$, $\epsilon$, $\epsilon$), where $\epsilon = \frac{\alpha + 2 \sigma}{4}$. 
Hence for a proper generalized curvature tensor, weak symmetry and Chaki pseudosymmetry are equivalent. We note that the solutions of a weakly $B$-symmetric manifold are not unique. We also note that if there is a solution of the form ($2\epsilon$, $\epsilon$, $\epsilon$, $\epsilon$, $\epsilon$), there may be another solution of the form ($\alpha$, $\beta$, $\beta$, $\gamma$, $\gamma$), $\beta \neq \gamma$. The solution of the form ($\alpha$, $\beta$, $\beta$, $\gamma$, $\gamma$) with different $\beta$ and $\gamma$ are studied by many authors (see, \cite{DB}, \cite{JANA3}, \cite{SK}).
We mention that for a proper generalized curvature tensor $B$ if the dimension of the space
$$\left\{\xi\in \chi^*(M):\xi(X_1)\, B(X_2,X_3)X_4+\xi(X_2)\, B(X_3,X_1)X_4+\xi(X_3)\, B(X_1,X_2)X_4=0\right\}$$
is zero, where $\chi^*(M)$ is Lie algebra of all smooth $1$-forms on $M$, then the solution of a weakly $B$-symmetric manifold 
is uniquely determined as ($2\epsilon$, $\epsilon$, $\epsilon$, $\epsilon$, $\epsilon$), 
where $\epsilon = \frac{\alpha + 2 \sigma}{4}$.\\
\indent We now discuss the results for a $(0,2)$-tensor $Z$ to be weakly symmetric. 
For this purpose, we need the following obvious properties of a $(0,2)$-tensor $Z$.
\begin{lem}\label{lem1}
Let $(M,g)$, $n \geqslant 3$, be a semi-Riemannian manifold  
admitting an $(0,2)$-tensor field $Z$ and an $1$-form $\theta$.\\
$(1)$ If $\theta(X_1)Z(X_2,X_3)+\theta(X_2)Z(X_1,X_3) = 0,$ then either $Z=0$ or $\theta =0$.\\
$(2)$ If $Z$ is a symmetric tensor such that $\theta(X_1)\, Z(X_2,X_3)+\theta(X_2)\, Z(X_3,X_1)+\theta(X_3)\, Z(X_1,X_2) = 0,$ then either $Z=0$ or $\theta =0$.\\
$(3)$ If $Z$ is a symmetric tensor such that $\theta(X_1)\, Z(X_2,X_3)-\theta(X_2)\, Z(X_1,X_3) = 0,$ then either rank$(Z)\leqslant 1$ or $\theta =0$.
\end{lem}
\indent Now we consider a weakly $Z$-symmetric semi-Riemannian manifold
$(M,g)$, $n \geqslant 3$, with solution ($\delta$, $\eta$, $\lambda$). Then
\be\label{eq3.2}
(\nabla_X Z)(X_1,X_2) = \delta(X)\, Z(X_1,X_2) + \eta(X_1)\, Z(X,X_2) + \lambda(X_2)\, Z(X_1,X),
\ee
Now changing the position of $X$, $X_1$ and $X_2$, and combining the resultant equations as require, from Lemma \ref{lem1} we have the following:
\begin{thm}\label{thm3.1}
Let $(M,g)$, $n \geqslant 3$, be a weakly $Z$-symmetric manifold with solution $(\delta$, $\eta$, $\lambda)$. Then we have\\
\indent$(1)$ if $Z$ is symmetric, then $(i)$ there exists a solution $(\delta$, $\nu$, $\nu)$ such that $\nu = \frac{\eta+\lambda}{2}$,\\
\indent\indent\hspace{1.57in} $(ii)$ $\eta = \lambda$ if rank$(Z)$ $>1$ $($\cite{mmid}, \cite{mmwsn}$)$.\\
\indent$(2)$ if $Z$ is non-zero and skew-symmetric, then $\eta = \lambda$.\\
\indent$(3)$ if $Z$ is symmetric and Codazzi type with rank$(Z)>1$ then $\delta = \eta = \lambda$.\\
\indent$(4)$ if $Z$ is symmetric and cyclic parallel then $\delta+\eta+\lambda =0$ and hence there exists a solution\\
\indent \ \ \ \ $(2 \zeta$, -$\zeta$, -$\zeta)$ such that $\zeta = \frac{\eta+\lambda}{2}$.
\end{thm}
\begin{rem}
Thus we conclude that the solution $(\delta$, $\eta$, $\lambda)$ of a $(WRS)_n$ turns into $(\delta$, $\nu$, $\nu)$ 
such that $\nu = \frac{\eta+\lambda}{2}$. The solution of the form $(\delta$, $\eta$, $\lambda)$ 
of a $(WRS)_n$ with different $\eta$ and $\lambda$ are studied by many authors (see, \cite{DEGS}, \cite{mmid}, \cite{SK} and also references therein).
\end{rem}
\begin{thm} $($\cite{ManSuh02}, \cite{ManSuh12b}$)$ 
If a $(WS)_n$ admits a solution other than of the form ($2\phi$, $\phi$, $\phi$, $\phi$, $\phi$) then the curvature 2-form $\Omega_{(R)l}^m$ is recurrent. Moreover the result is true for any generalized curvature tensor.
\end{thm}
\noindent \textbf{Proof}: It can be easily shown that a manifold satisfying \eqref{eq3.1} also satisfies \eqref{b3}  for $T=B$ with the corresponding 1-form $\alpha - 2\beta \neq 0$ or $\alpha - 2\gamma \neq 0$ (from hypothesis) and thus the curvature 2-form $\Omega^m_{(B)l}$ is recurrent.
\begin{thm} $($\cite{ManSuh12b}$)$ 
In a $(WRS)_n$ the Ricci 1-form $\Omega_{(S)l}$ is recurrent. Moreover the result is true for any symmetric $(0,2)$ tensor $Z$.
\end{thm}
\noindent \textbf{Proof}: A manifold satisfying \eqref{wrsntb}, satisfies \eqref{b4} if $Z$ is symmetric. Thus the 1-form $\Lambda_{(Z)l}$ is recurrent (see \cite{ManSuh12b} Remark 2.6.).
%
\begin{defi}\label{newdefiI}
A semi-Riemannian manifold $(M,g)$, $n \geqslant 3$, 
admitting an $(0,k)$-tensor field $T$, $k\ge 2$, 
is said to be weakly $T$-symmetric of type-$I$ if
\bea\label{newdefi1}
(\nabla_{X_1}T)(X_2,X_3, \ldots ,X_{k+1}) = \sum_p \stackrel{\stackrel{p}{}}{\alpha}(X_{p(1)})\, T(X_{p(2)},X_{p(3)}, \ldots ,X_{p(k+1)}),
\eea
where $\stackrel{\stackrel{p}{}}{\alpha}$ are associated $1$-forms and the sum includes all permutation $p$ 
over the set $(1,2, \ldots ,k+1)$.
\end{defi}
We note that this defining condition of weakly $T$-symmetric manifold is due to Ewert-Krzemieniewski \cite{sti}.
\begin{defi}\label{newdefiII}
A semi-Riemannian manifold $(M,g)$, $n \geqslant 3$, 
admitting an $(0,k)$-tensor field $T$, $k\ge 2$, 
is said to be
weakly $T$-symmetric of type-$I\!I$ if
\bea\label{newdefi2}
(\nabla_X T)(X_1,X_2, \ldots ,X_k) &=& \alpha(X)\, T(X_1,X_2, \ldots ,X_k)\\\nonumber
&+& \sum^k_{i=1}\pi_i(X_i)\, T(X_1,X_2, \ldots ,\underset{i-th\ place}{X}, \ldots ,X_k),
\eea
where $\alpha$ and $\pi_i$, $i=1,2, \ldots ,k$ are associated $1$-forms.
\end{defi}
We note that this defining condition of weakly $T$-symmetric manifold is due to Tam$\acute{\mbox{a}}$ssy and Binh \cite{tb}.
\begin{defi}\label{newdefiIII}
A semi-Riemannian manifold $(M,g)$, $n \geqslant 3$, 
admitting an $(0,k)$-tensor field $T$, $k\ge 2$, 
is said to be
weakly $T$-symmetric of type-$I\!I\!I$ if
\bea\label{newdefi3}
\nabla_X T = \alpha\otimes T  - \pi_{_X} \cdot T,
\eea
where $\alpha$ and $\pi$ are called associated $1$-forms.
\end{defi}
Now if we set $\alpha = 2 \pi$ in $(\ref{newdefi3})$ then a weakly $T$-symmetric manifold of type-$I\!I\!I$ takes the form 
of Chaki $T$-pseudosymmetric manifold. Thus the defining condition of a Chaki $T$-pseudosymmetric manifold for a $(0,k)$-tensor is given by
\bea\label{newpsn}
\nabla_X T = 2 \pi\otimes T  - \pi_{_X} \cdot T.
\eea
again if we set $\pi =0$ in $(\ref{newdefi3})$ then a weakly $T$-symmetric manifold of type-$I\!I\!I$ takes the form 
of $T$-recurrent manifold. Thus the defining condition of a $T$-recurrent manifold for a $(0,k)$-tensor is given by
\bea\label{newkn}
\nabla_X T = 2 \alpha\otimes T.
\eea
We note that weak symmetry of type-$I\!I\!I$ is a special case of type-$I\!I$ and type-$I\!I$ is a special case of type-$I$. 
Moreover, if we consider $k=4$ in the above definition of weakly $T$-symmetric of type-$I\!I$ and Chaki pseudosymmetric manifold we get the equations \eqref{wsntb} and \eqref{cpsn} respectively.
\begin{thm}\label{th3.2}
Let $T$ be an $(0,k)$-tensor, $k\geqslant 2$, on a semi-Riemannian manifold $(M,g)$, $n \geqslant 3$,
skew symmetric in $i$-th and $j$-th indices $\big{(}i,j\in \{1,2, \ldots , k\}, \ i\neq j\big{)}$. 
If $(M,g)$ is a weakly $T$-symmetric manifold of type-$I\!I$ then $\pi_i = \pi_j$.
\end{thm}
\noindent\textbf{Proof:} Proof is similar to the first part of the proof of Theorem 1 of \cite{pra} and hence we omit it.
\begin{thm}\label{th3.3}
Let $T$ be an $(0,k)$-tensor, $k\geqslant 2$, on a semi-Riemannian manifold $(M,g)$, $n \geqslant 3$,
symmetric with respect to $p$-number $(2p\leqslant k)$ 
of indices with other $p$-number of indices taken together. 
If $(M,g)$ is a weakly $T$-symmetric manifold of type-$I\!I$ then
each corresponding $1$-forms can be replaced by same in pair.
\end{thm}
\noindent\textbf{Proof:} Interchanging the $p$-number of indices with corresponding $p$-number 
of indices taken together in (\ref{newdefi2}) and adding the resulting equation with (\ref{newdefi2}), 
and then using the symmetry of $T$, we obtain the result.
%
\section{\bf{Condition of weak symmetry and Chaki pseudosymmetry to be Deszcz pseudosymmetry}}\label{Main}
%
In this section we deduce the condition for a Chaki pseudosymmetry and weak symmetry of type-$I\!I\!I$ 
to be Deszcz pseudosymmetry for a $(0,k)$-tensor $T$, $k\ge 2$.\\
Let $(M,g)$, $n \geqslant 3$, be a weakly $T$-symmetric of type-$I\!I\!I$ semi-Riemannian manifold
with defining condition (\ref{newdefi3}). Then by a straightforward calculation we have
\bea
&&(R \cdot T)(X_1,X_2,\ldots,X_k;X,Y)\\\nonumber
&&= 2d\alpha(X,Y)\otimes T(X_1,X_2,\ldots,X_k) + Q(J,T)(X_1,X_2,\ldots,X_k;X,Y),
\eea
\noindent where $J = \pi\otimes \pi  - \nabla \pi$ and $d\alpha$ denotes the exterior derivative of $\alpha$. 
This leads to the following:
\begin{thm}\label{main3}
A weakly $T$-symmetric semi-Riemannian manifold $(M,g)$, $n \geqslant 3$, of type-$I\!I\!I$ is\\
$(i)$ $T$-semisymmetric if and only if $2 d\alpha\otimes T + Q(J,T)$ is zero,\\
$(ii)$ Deszcz $T$-pseudosymmetric $($resp., Ricci generalized $T$-pseudosymmetric$)$ if and only if $2 d\alpha\otimes T + Q(J,T)$ is linearly dependent with $Q(g,T)$ $($resp., $Q(S,T))$.
\end{thm}
If $T=R$ then we get the results for a $(WS)_n$ and for $T= S$ we get the results for a $(WRS)_n$. 
Now as a direct consequence of the above theorem we can state the following:
\begin{cor}\label{cor4.1}
A weakly $T$-symmetric semi-Riemannian manifold $(M,g)$, $n \geqslant 3$, of type-$I\!I\!I$ is\\
$(i)$ $T$-semisymmetric if the associated $1$-form $\alpha$ is closed and $(\pi\otimes \pi - \nabla \pi) = 0$,\\
$(ii)$ Deszcz $T$-pseudosymmetric $($resp., Ricci generalized $T$-pseudosymmetric$)$ if the associated $1$-form $\alpha$ is closed and $(\pi\otimes \pi - \nabla \pi)$ is proportional to $g$ $($resp., $S)$,\\
$(iii)$ $T$-semisymmetric if the associated $1$-form $\alpha$ is closed and $Q(J,T)$ is zero,\\
$(iv)$ Deszcz $T$-pseudosymmetric $($resp., Ricci generalized $T$-pseudosymmetric$)$ if the associated $1$-form $\alpha$ is closed and $Q(J,T)$ is linearly dependent with $Q(g,T)$ $($resp., $Q(S,T))$.
\end{cor}
\begin{cor}\label{cor4.3}
A weakly $B$-symmetric semi-Riemannian manifold $(M,g)$, $n \geqslant 3$, of type-$I\!I\!I$ is\\
$(i)$ $B$-semisymmetric if the associated $1$-form $\alpha$ 
is closed and $\pi$ is such that $B = L_1 \, J \wedge J$, where $J = \pi\otimes \pi - \nabla \pi$,\\
$(ii)$ Deszcz $B$-pseudosymmetric if the associated $1$-form $\alpha$ is closed 
and $\pi$ is such that $B = L_1\, D_1\wedge D_1$, where $D_1 = L_2 g - J$ $($in this case 
$R\cdot  B = L_2\, Q(g,B))$,\\
$(iii)$ Ricci generalized $B$-pseudosymmetric if the associated $1$-form $\alpha$ 
is closed and $\pi$ is such that $B = L_3\, D_2\wedge D_2$, where $D_2 = L_4 S - J$ $($in this case $R\cdot B = L_4\, Q(S,B))$,\\
where $L_1, L_2, L_3$ and $L_4$ are smooth functions chosen suitably.
\end{cor}
Since Chaki $T$-pseudosymmetry is a special case of weak $T$-symmetry, from Theorem \ref{main3} 
we can state the results for a manifold with Chaki pseudosymmetry to be Deszcz pseudosymmetric as follows:
\begin{cor}\label{cor4.4}
A Chaki $T$-pseudosymmetric semi-Riemannian manifold $(M,g)$, $n \geqslant 3$,
with associated $1$-form $\phi$, is\\
$(i)$ $T$-semisymmetric if and only if $4 d\phi\otimes T + Q(\phi\otimes \phi - \nabla \phi,T)$ is zero,\\
$(ii)$ Deszcz $T$-pseudosymmetric $($resp., Ricci generalized $T$-pseudosymmetric$)$ if and only if $4 d\phi\otimes T + Q(\phi\otimes \phi - \nabla \phi,T)$ is linearly dependent with $Q(g,T)$ $($resp., $Q(S,T))$.
\end{cor}
\begin{cor}\cite{DES}\label{cor4.5}
A Chaki $T$-pseudosymmetric semi-Riemannian manifold $(M,g)$, $n \geqslant 3$,
with associated $1$-form $\phi$, is\\
$(i)$ $T$-semisymmetric if $\phi$ is closed and $(\phi\otimes \phi - \nabla \phi)=0$,\\
$(ii)$ Deszcz $T$-pseudosymmetric $($resp., Ricci generalized $T$-pseudosymmetric$)$ if $\phi$ is closed 
and the tensor $(\phi\otimes \phi - \nabla \phi)$ is proportional to $g$ $($resp., $S)$.
\end{cor}
We note that if at every point $x\in M$ the tensor ($\phi\otimes \phi - \nabla \phi$) 
is zero or proportional to the metric tensor $g$ or the Ricci tensor $S$, then closedness of $\phi$ is obvious. 
Thus the condition of closedness of $\phi$ in the above Corollary is not required. 
We note that the conditions of Corollary \ref{cor4.5} are not necessary (see, Example 5.4 of the last section).
We also mention that 
\cite{DES} (Section 5.2) contains some comments related to pseudosymmetric (resp. Ricci-pseudosymmetric) manifolds
and Chaki pseudosymmetric (resp. pseudo-Ricci symmetric) manifolds. 
\begin{cor}\label{cor4.6}
A Chaki $T$-pseudosymmetric 
semi-Riemannian manifold $(M,g)$, $n \geqslant 3$, 
with associated $1$-form $\phi$, is\\
$(i)$ $T$-semisymmetric if $\phi$ is closed and $Q(\phi\otimes \phi - \nabla \phi,T)$ is zero,\\
$(ii)$ Deszcz $T$-pseudosymmetric $($resp., Ricci generalized $T$-pseudosymmetric$)$ if $\phi$ is closed and $Q(\phi\otimes \phi - \nabla \phi,T)$ is linearly dependent with $Q(g,T)$ $($resp., $Q(S,T))$.
\end{cor}
\begin{cor}\label{cor4.7}
A Chaki $B$-pseudosymmetric 
semi-Riemannian manifold $(M,g)$, $n \geqslant 3$,
with associated $1$-form $\phi$, is\\
$(i)$ $B$-semisymmetric if $\phi$ is such that $B = L_1\, H \wedge H$, 
where $H = \phi\otimes \phi - \nabla \phi$,\\
$(ii)$ Deszcz $B$-pseudosymmetric if $\phi$ is such that $B = L_1\, D_1\wedge D_1$, 
where $D_1 = L_2 g - H$ $($in this case $R\cdot B = L_2\, Q(g,B))$,\\
$(iii)$ Ricci generalized $B$-pseudosymmetric if $\phi$ is such that $B = L_3\, D_2\wedge D_2$, 
where $D_2 = L_4 S - H$ $($in this case $R \cdot B = L_4\, Q(S,B))$,\\
where $L_1, L_2, L_3$ and $L_4$ are smooth functions chosen suitably.
\end{cor}
We mention that a recurrent manifold \cite{Ag} is a special form of weakly symmetric manifold with $\beta=\gamma=0$. 
Now if $T$ is a $(0,k)$-tensor then in any $T$-recurrent Riemannian manifold for the function 
$f = T_{i_1 i_2 \ldots i_k}T^{i_1 i_2 \ldots i_k}$ we have $f_{,l} = 2 f \alpha_l $ and as a consequence
$ \alpha_l = \frac{1}{2} (\log |f|) _{, l}$ on the open subset of all points of $M$ at which $f \neq 0$. 
This implies that the associated $1$-form $\alpha$ is locally gradient and hence closed. 
Then from Theorem \ref{main3} we get the following:
\begin{cor}\label{cor4.8}
Every $T$-recurrent Riemannian manifold $(M,g)$, $n \geqslant 3$, is $T$-semisymmetric.
\end{cor}
\indent We note that if $T$ is any generalized curvature tensor then the corresponding results of the above corollary are reported in many papers.
%
\begin{defi} $($\cite{MIK}, \cite{sou}, \cite{shi}, \cite{yano1}, \cite{yano2}$)$ 
A vector field $V$ on a 
semi-Riemannian manifold $(M,g)$, $n \geqslant 3$, 
is said to be torseforming 
if it satisfies the equation of the form $\nabla_X V = a X + \tau(X) V$, $\forall X \in \chi(M)$, 
where $a$ is a scalar and $\tau$ is an $1$-form. If $\omega$ is the corresponding $1$-form of $V$, i.e., $\omega(X) = g(X,V)$ then
$$(\nabla_X \omega)(Y) = a g(X,Y) + \tau(X) \omega(Y), \forall X, Y \in \chi(M).$$
\end{defi}
\noindent The torseforming vector field $V$ on a 
semi-Riemannian manifold $(M,g)$, $n \geqslant 3$, 
is called \cite{MIK}\\
\indent(1) recurrent, if $a = 0$,\\
\indent(2) concircular, if $\tau$ is a gradient $1$-form, i.e., $\tau = dh$, $h$ being a scalar,\\
\indent(3) convergent, if it is concircular, and $a$ is a constant multiple of $e^{h}$,\\
\indent(4) proper concircular \cite{yano1}, if $\tau$ is closed.\\
We now state some fundamental well known \cite{shi} results on torseforming vector fields.\\
\indent (i) A non-recurrent torseforming vector field $V$ is non-isotropic i.e. $g(V,V) \neq 0$.\\
\indent (ii) The constant multiple of a torseforming vector field is a torseforming vector field. 
However, if $V$ is a torseforming vector field, then $f V$ is not necessarily a torseforming vector field for any smooth function $f$.\\
\indent (iii) Any unit torseforming vector field $\hat{V}$ 
is of the form $\nabla_X \hat{V} = a \left( X + \hat{\omega}(X) \hat{V} \right)$, $\forall X \in \chi(M)$, 
where $\hat{\omega}$ is the $1$-form corresponding to the unit vector field $\hat{V}$, 
here `$a$' is the associated scalar of the unit torseforming vector field.\\
\indent (iv) If $V$ is a non-recurrent torseforming vector field such that $da$ and $\omega$ are linearly dependent 
then $\omega$ is closed if and only if $\tau$ is closed. 
Also if $\omega$ and $\tau$ both are closed then $da$ and $\omega$ are linearly dependent.\\
\indent (v) If $g(V,V)$ is a non-zero constant, then $\omega$ is closed.\\
\indent (vi) If $\omega$ is closed, then $\omega$ and $\tau$ are linearly dependent 
and hence $\nabla_X V = a\, X + b\, \omega(X) V$, consequently $(\nabla_X \omega)(Y) = a\, g(X,Y) + b\, \omega(X) \omega(Y)$, where $a$ and $b$ are the associated scalars of the closed torseforming vector field $V$.\\
Then the geometric significance of Corollary \ref{cor4.1} and \ref{cor4.5} are given by the following:
\begin{cor}\label{cor4.9}
If in a weakly $T$-symmetric 
semi-Riemannian manifold $(M,g)$, $n \geqslant 3$, 
of type-$I\!I\!I$ with associated $1$-forms $\alpha$ and $\pi$ 
are closed and the vector field corresponding to $\pi$ is torseforming with $b=1,$ 
then $(M,g)$ 
is a Deszcz $T$-pseudosymmetric.
\end{cor}
\begin{cor}\label{cor4.10}
If in a Chaki $T$-pseudosymmetric 
semi-Riemannian manifold $(M,g)$, $n \geqslant 3$, 
the associated $1$-form $\phi$ is closed and the corresponding vector field is torseforming with $b=1$, 
then $(M,g)$ is a Deszcz $T$-pseudosymmetric.
\end{cor}
Now we suppose that a generalized curvature tensor $B$ satisfies the condition
\bea\label{eq4.2}
(R(X_1,X_2) \cdot B)(X_3,X_4,X_5,X_6) & + & (R(X_3,X_4) \cdot B)(X_5,X_6,X_1,X_2)\\\nonumber
&+&(R(X_5,X_6)\cdot B)(X_1,X_2,X_3,X_4)=0.
\eea
We note that if we take $B=R$ then (\ref{eq4.2}) turns into Walker identity \cite{Ag}. 
By virtue of Walker identity in \cite{sti} (Theorem 1) it is shown that in a $(CPS)_n$ 
(and hence in a $(WS)_n$) the associated $1$-form (reduced 1-form) is closed. 
Hence if a Chaki $B$-pseudosymmetric manifold satisfies (\ref{eq4.2}), then its associated $1$-form is closed. 
Similarly it is easy to check that if the manifold is weakly $B$-symmetric 
of type-$I\!I\!I$ such that (\ref{eq4.2}) holds then the closedness of one associated $1$-form implies the closedness of the other. Thus for a weakly $B$-symmetric manifold $(M,g)$ with solution $(\alpha, \pi, \pi, \pi, \pi)$ 
and $B$ satisfying (\ref{eq4.2}), $(M,g)$ is Deszcz pseudosymmetric if any one of $\alpha$ 
or $\pi$ is closed and the corresponding vector field of $\pi$ is torseforming with the associated scalar $b= 1$. 
Also for a Chaki $B$-pseudosymmetric manifold $(M,g)$ with solution $(2 \pi, \pi, \pi, \pi, \pi)$ and $B$ satisfying (\ref{eq4.2}), $(M,g)$ is Deszcz pseudosymmetric if the corresponding vector field of $\pi$ is torseforming with the associated scalar $b= 1$. Again in a Chaki $T$-pseudosymmetric manifold $(M,g)$ if the corresponding vector field of the associated $1$-form 
is a unit proper torseforming with $a = 1$ as the associated scalar, then $(M,g)$ is a Deszcz $T$-pseudosymmetric. And in a weakly $T$-symmetric manifold $(M,g)$ if the associated $1$-form $\alpha$ is closed and the corresponding vector field of $\pi$ is unit torseforming with $a = 1$ as the associated scalar, then $(M,g)$ is a Deszcz $T$-pseudosymmetric. Hence we conclude that the closedness of the associated $1$-forms in Corollary \ref{cor4.9} and \ref{cor4.10} can be replaced by the condition (\ref{eq4.2}).\\
\indent A semi-Riemannian manifold $(M,g)$, $n \geqslant 4$, is said to be a manifold with pseudosymmetric Weyl tensor 
(\cite{DES}, \cite{DESZ})  
if the tensors $C \cdot C$ and $Q(g,C)$ are linearly dependent 
at every point of $M$. This is equivalent to  
\begin{eqnarray*}
C \cdot C = L_{C}\, Q(g,C)
\end{eqnarray*}
on $U_{C}$, where $L_{C}$ is some function on this set.             
\section{\bf{Some examples of Deszcz and Chaki pseudosymmetric manifolds}}\label{Example}
\noindent\textbf{Example 5.1.} 
Let $M$ be a non-empty open connected subset of $\mathbb{R}^5$ endowed with the metric $g$ defined by  
\begin{eqnarray*}
g_{ij}dx^{i} dx^{j} &=& e^{x^1}(dx^1)^2 + e^{x^1} \left( e^{x^5}(dx^2)^2 + (dx^3)^2 + (dx^4)^2+ (dx^5)^2 \right),\ \ i, j = 1,2, \ldots ,5.
\end{eqnarray*}
Then by a straightforward calculation we can easily evaluate the components of Riemann-Christoffel curvature tensor, Ricci tensor 
and the covariant derivative of the curvature tensor of $(M,g)$. Again using these values we can evaluate the non-zero components of $R\cdot R$, $Q(g,R)$ and $Q(S,R)$ easily. Let us now consider the $1$-form $\phi$ as follows: 
\begin{eqnarray*}
\phi_{i}(x)=\left\{\begin{array}{clc}
&-\frac{1}{2}&\ \ \ \ \mbox{for}\ \ \ i=1\\
&0&\ \ \ \ \mbox{otherwise.}
\end{array}\right.
\end{eqnarray*}
Then for this $1$-form $\phi$, 
the manifold $(M,g)$ satisfies \eqref{cpsn} and thus becomes a Chaki pseudosymmetric manifold. But this is neither Deszcz pseudosymmetric nor Ricci generalized pseudosymmetric. Moreover, the following relations also are satisfied on $(M,g)$:  
$\kappa  = \frac{7}{2} e^{ - x^{1}}$ and
\begin{eqnarray*}
R &=& 
\frac{65}{36} e^{x^{1}}  S\wedge S 
- \frac{34}{9} e^{2 x^{1}} S\wedge S^{2} 
+ \frac{20}{9} e^{3 x^{1}} S^{2}\wedge S^{2},\\
S^{2}\wedge S^{2} 
&=& 
\frac{7}{2}  e^{ - x^{1}} S\wedge S^{2} 
- \frac{49}{16}  e^{ - 2 x^{1}} S\wedge S 
- \frac{3}{2} e^{ - 2 x^{1}} g\wedge S^{2} 
+ \frac{21}{8}  e^{ - 3 x^{1}} g\wedge S 
-  \frac{9}{8} e^{ - 4 x^{1}} G ,\\
C\cdot C 
&=& 
- \frac{1}{24} e^{ - x^{1}} Q(g,C),\\ 
R\cdot R - Q(S,R) 
&=& 
\frac{25}{9} e^{2 x^{1}} Q(S,S\wedge S^{2}) 
+ \frac{28}{9} e^{3 x^{1}} Q(S^{2},S\wedge S^{2}).
\end{eqnarray*}
We recall that the Ricci operator $\mathcal{S}$ and the $(0,2)$-tensor $S^{2}$ 
of a semi-Riemannian manifold $(M,g)$ are defined by
$g(\mathcal{S} X,Y) = S(X,Y)$ and   
$S^{2}(X,Y) = S( \mathcal{S}X,Y)$, respectively, where $X,Y \in \chi(M)$. 
We mention that in \cite{kun} (Example 4.1) an example of 
a spacetime with the Riemann-Christoffel curvature tensor $R$ expressed 
by a linear combination of the tensors $S\wedge S$, $S\wedge S^{2}$ and $S^{2}\wedge S^{2}$
is given.
\begin{defi}
A semi-Riemannian manifold $(M,g)$, $n \geqslant 4$, is said to be Roter type \cite{DA} if its curvature tensor $R$ is expressed as the linear combination of $g\wedge g$, $g\wedge S$ and $S\wedge S$ i.e.,
\be\label{eqr}
R=N_1\, g\wedge g + N_2\, g\wedge S + N_3\, S\wedge S,
\ee
where $N_1$, $N_2$ and $N_3$ are some smooth functions on $M$.
\end{defi}
Evidently, on every conformally flat semi-Riemannian manifold $(M,g)$, $n \geqslant 4$, we have
\begin{eqnarray*}
R &=&\frac{1}{n-2} g \wedge S - \frac{\kappa}{2(n-1)(n-2)} g \wedge g .
\end{eqnarray*}
Thus (\ref{eqr}) is satisfied.
\begin{defi}
A semi-Riemannian manifold $(M,g)$, $n \geqslant 4$, is said to be \textit{generalized Roter type} 
if its curvature tensor $R$ is expressed as the linear combination of 
$S\wedge S$, $S\wedge S^2$, $g\wedge S$, $g\wedge S^2$, $g\wedge g$ and $S^2\wedge S^2$ i.e.,
\be\label{eqgr}
R= L_1\, S\wedge S + L_2\, S\wedge S^2 + L_3\, g\wedge S + L_4\, g\wedge S^2 + L_5\, g\wedge g + L_6\, S^2\wedge S^2,
\ee
where $L_i$, $1, 2, \ldots , 6$ are some smooth functions on $M$.
\end{defi}
We note that any Roter type manifold is a generalized Roter type but not conversely. 
We mention that non-Roter type manifolds with the curvature tensor having  
a decomposition of the form (\ref{eqgr}) were already investigated in \cite{Saw3} 
and very recently in \cite{DGJP-TZ}, \cite{kun}. 
Namely, \cite{Saw3} contains results on hypersurfaces in space forms having
curvature tensors of the form (\ref{eqgr}). 
As it was shown in Section 2 of \cite{DGJP-TZ},
the $4$-dimensional manifold
presented in Section 4 of \cite{DebKon} 
has the curvature tensor of the form (\ref{eqgr}). 
Some spacetimes also satisfy (\ref{eqgr}) (\cite{kun}, Example 4.1).\\
\indent It is easy to check that for the aforesaid manifold (\ref{eqr}) does not hold but (\ref{eqgr}) holds with
\beb
&L_1 = -\frac{1}{16} e^{-2 x^1} \left(80 e^{3 x^1}-49 L_6\right), \ \ 
L_2 = \frac{1}{2} e^{-x^1} \left(8 e^{3 x^1}-7 L_6\right),&\\
&\frac{24}{7} e^{3 x^1} L_3 = -6 e^{2 x^1} L_4 = -16 e^{4 x^1} L_5 = \left(20 e^{3 x^1}-9 L_6\right),&
\eeb
and arbitrary smooth function $L_6$.
Hence the manifold $(M,g)$ under consideration is generalized Roter type but not Roter type. 
We note that $(M,g)$ do not satisfy (\ref{b1}), (\ref{b2}) and (\ref{b3}) for $T = C, P, K$ and $conh(R)$, 
and also do not satisfy (\ref{b2}) for $T=R$ but (\ref{b4}) holds for $Z = S$ 
with $\alpha = \left\{-\frac{1}{2},0,0,0,0\right\}$. 
Hence the manifold is of recurrent Ricci $1$-form but not of recurrent curvature $2$-form for $R, C, P, K$ and $conh(R)$.\\
\noindent\textbf{Example 5.2.} 
Let $M$ be a non-empty open connected subset of $\mathbb{R}^4$, where $x^1 > 0$, endowed with the metric 
$g$ defined by 
\begin{eqnarray*} 
g_{ij}dx^{i} dx^{j} &=& x^1 \left( (dx^1)^2 + (dx^2)^2 + (dx^3)^2 + (dx^4)^2 \right), \ \ i, j = 1,2,3,4.
\end{eqnarray*}
Evidently, $(M,g)$ is a conformally flat manifold with $\kappa  = -\frac{3}{2 (x^{1})^3} $. Thus (\ref{eqr}) is satisfied. The manifold $(M,g)$ is not a Chaki pseudosymmetric manifold. However, this is a Deszcz pseudosymmetric and Ricci generalized pseudosymmetric manifold with
$$R \cdot R = -\frac{1}{2 (x^1)^3} Q(g,R) = Q(S,R).$$
Moreover, in the above manifold $S\wedge S=0 $ , $S\wedge S^2=0 $ and $S^2\wedge S^2=0$. We note that for $T=R$, $(M,g)$ fulfills (\ref{b1}) but does not fulfill (\ref{b2}). 
Again for $T= P, K$ and $conh(R)$, $(M,g)$ do not realize (\ref{b1}) and (\ref{b2}) but realizes (\ref{b3}) 
for $\alpha=\left(-\frac{1}{x^1},0,0,0\right)$, $\left(-\frac{1}{x^1},0,0,0\right)$ 
and $\left(-\frac{3}{x^1},0,0,0\right)$ respectively. 
We also note that $M$ does not satisfy (\ref{b4}) for $Z = S$. 
Hence the manifold is of recurrent curvature $2$-form for $P, K$ and $conh(R)$ 
but the curvature $2$-form corresponding to $R$ and the Ricci $1$-form are not recurrent.\\
\noindent\textbf{Example 5.3.}  
We define on $\mathbb R^4$ the metric $g$ by
\begin{eqnarray*}
g_{ij}dx^i dx^j &=& a^2 \left( -(dx^1)^2 + \frac{1}{2}e^{2x^1}(dx^2)^2 - (dx^3)^2 + (dx^4)^2 + 2e^{x^1}dx^2 dx^4 \right),\ \ i, j = 1,2,3,4,
\end{eqnarray*}
where $a$ is a non-zero constant. 
The manifold $(\mathbb R^4,g)$ is called the G\"{o}del spacetime \cite{KG}.
It is well-known that the G\"{o}del spacetime is a non-conformally flat manifold with
the Ricci tensor $S$ of rank one. Thus the G\"{o}del spacetime is a quasi-Einstein manifold.
The G\"{o}del spacetime is a manifold with  pseudosymmetric Weyl tensor (\cite{kun}): $C \cdot C = \frac{\kappa }{6} Q(g,C)$.
Moreover,  it can be seen that this manifold is neither a Chaki pseudosymmetric manifold 
nor Deszcz pseudosymmetric but Ricci generalized pseudosymmetric \cite{kun}. 
For more details about the
other curvature properties of G\"{o}del metric see \cite{kun}. 
We note that G\"{o}del spacetime do not satisfy (\ref{b1}), (\ref{b2}) and (\ref{b3}) 
for $T= C, P, K, conh(R)$ except (\ref{b1}) for $T=K$, and also (\ref{b2}), (\ref{b4}) 
does not hold for $T=R$, $Z = S$, respectively. So G\"{o}del spacetime 
is neither of recurrent curvature $2$-form for $R, C, P, K, conh(R)$ nor of recurrent Ricci $1$-form.\\
\noindent\textbf{Example 5.4.} 
Let $M$ be a non-empty open connected subset of $\mathbb{R}^4$ endowed with the metric 
$g$ defined by  
\begin{eqnarray*}
g_{ij}dx^{i} dx^{j} &=& (e^{x^1}+1)(dx^1)^2 + e^{x^1} \left( (dx^2)^2 + (dx^3)^2 + (dx^4)^2 \right),\ \ i,j = 1,2,3,4.
\end{eqnarray*}
Evidently, $(M,g)$ is a conformally flat manifold with scalar curvature $\kappa=\frac{3(2+e^{x^1})}{2(1+e^{x^1})^2}$. Let us now consider the $1$-form $\phi$ as follows: 
\begin{eqnarray*}
\phi_{i}(x)=\left\{\begin{array}{clc}
&-\frac{e^{x^1}}{2 \left(e^{x^1}+1\right)}&\ \ \ \ \mbox{for}\ \ \ i=1\\
&0&\ \ \ \ \mbox{otherwise.}
\end{array}\right.
\end{eqnarray*}
Then for this $1$-form $\phi$, the manifold $(M,g)$ is a Chaki pseudosymmetric manifold. 
Now $H = \phi\otimes \phi -\nabla \phi$ is given by
$$H_{11} = \frac{e^{x^1}}{2(1+e^{x^1})^2}, \ \ \ H_{22}= H_{33} = H_{44} = \frac{e^{2 x^1}}{4(1+e^{x^1})^2}.$$
Then it is clear that $H$ is not proportional to $g$ or $S$. Thus 
the sufficient
condition 
for Corollary \ref{cor4.5} does not hold, so now we can not get any conclusion for this manifold to be $(DPS)_4$ 
or Ricci generalized pseudosymmetric. Again
$$R = \frac{2 (e^{x^1}+1)^3}{(e^{x^1}-1)^2}\, D_1\wedge D_1 = \frac{2 (e^{x^1}+1)^3}{(3+e^{x^1})^2}\, D_2\wedge D_2,$$
where $D_1 = \frac{1}{4(1+e^{x^1})^2}g - H$ and $D_2 = S - H$. Thus by Corollary \ref{cor4.7} we have
$$R\cdot R = \frac{1}{4(1+e^{x^1})^2} Q(g,R) = Q(S,R)$$
i.e., the manifold $(M,g)$ is $(DPS)_4$ and also Ricci generalized pseudosymmetric. 
This entails that the conditions in Corollary \ref{cor4.5} are not necessary 
for a $(CPS)_4$ to be either a $(DPS)_4$ or a Ricci generalized pseudosymmetric manifold. 
It is easy to check that all the results of Theorem \ref{main3} hold for this example. 
Moreover, for the above manifold $(M,g)$ the curvature tensor can be expressed as  
$$R = N_1\, g\wedge g + N_2\, g\wedge S + N_3\, S\wedge S,$$
where 
$N_1 = -\frac{2+e^{x^1}}{4(1+e^{x^1})^2} + \frac{(3+2e^{x^1})^2}{16(1+e^{x^1})^4}N_3, \ \ \ \ 
N_2 = \frac{1}{2} + \frac{(3+2e^{x^1})}{4(1+e^{x^1})^2}N_3.$
Thus the manifold $(M,g)$ is Roter type and hence generalized Roter type. Again we note that for $T=R$, $(M,g)$ satisfies (\ref{b1}) and does not satisfy (\ref{b2}). For $T = P, K$ and $conh(R)$, $(M,g)$ do not satisfy (\ref{b1}), (\ref{b2}) but satisfies (\ref{b3}) for $\alpha=\left(-\frac{3+e^{x^1}}{1+e^{x^1}},0,0,0\right)$, $\left(-\frac{3+e^{x^1}}{1+e^{x^1}},0,0,0\right)$ and $\left(-\frac{e^{x^1}(3+e^{x^1})}{2+3e^{x^1}+e^{2x^1}},0,0,0\right)$ respectively. It also realizes (\ref{b4}) for $Z= S$ with $\alpha = \left(-\frac{e^{x^1} \left(e^{x^1}+3\right)}{5 e^{x^1}+2 e^{2 x^1}+3},0,0,0\right)$. Hence the manifold is of recurrent curvature $2$-form for $P, K$ and $conh(R)$ but not for $R$. Also it is of recurrent Ricci $1$-form.\\
\noindent\textbf{Example 5.5.} 
Let $M$ be an open connected subset of $\mathbb{R}^5$ endowed with the metric 
$g$ of the form  (Theorem 2.1; \cite{Kowalski01})
\begin{eqnarray*}
dx^{2} + dy^{2} + du^{2} + dv^{2} + \rho ^{2}  (x du - y dv + dz )^{2} ,
\end{eqnarray*}
where $\rho $ is a non-zero constant. $(M,g)$ is a non-conformally flat manifold. 
Its Ricci tensor $S$ is not proportional to $g$, cyclic parallel 
and 
$(S - \frac{\kappa }{2} g )\wedge ( S - \frac{\kappa }{2} g )  = 0$, $\kappa = \rho ^{2}$,
hold on $M$. Thus $(M,g)$ is a quasi-Einstein manifold
(cf. \cite{GLOG}, Lemma 3.1) with
$S=\alpha\, g + \beta\, \eta\otimes \eta$,
where $\alpha = \frac{\kappa }{2}$, $\beta = -\frac{3\kappa }{2}$ and $\eta = ( 0, 0, -\rho, -x\rho, y\rho )$.
This manifold is not Chaki pseudosymmetric 
but we have some pseudosymmetric type conditions on $M$: $R \cdot R =  - \frac{\kappa }{4}\, Q(g,R)$,
$C \cdot R =  - \frac{1}{3}\, Q(S,C) - \frac{\kappa }{3}\, Q(g,C)$,
$C \cdot S = 0$,
$C \cdot C = C \cdot R$, 
$R \cdot C - C \cdot R = \frac{1}{3}\, Q(S,R) + \frac{\kappa }{12}\, Q(g,R)$,
$P \cdot R$ and $Q(g,R)$ are linearly independent but $P\cdot S = - \frac{\kappa }{4}\, Q(g,S)$,
$K\cdot R = - \frac{3\kappa }{10}\, Q(g,R)$, 
$conh(R) \cdot R$ and $Q(g,R)$ are linearly independent but $conh(R) \cdot S = - \frac{\kappa }{12}\, Q(g,S)$ and in more general
\beb
\left(\alpha_2-\frac{4\alpha_1}{\kappa}\right) R\cdot R &+&
\left(3\alpha_2 + 4\alpha_4 - \frac{4 \alpha_3}{\kappa}\right) R\cdot C - 
\left(\beta + 3\alpha_2 + 3\alpha_4\right) C\cdot R + \beta\, C\cdot C\\
&=& \alpha_1 \, Q(g, R) + \alpha_2 \, Q(S, R) + \alpha_3 \, Q(g, C) + \alpha_4 \, Q(S, C),
\eeb
where $\alpha_1$, $\alpha_2$, $\alpha_3$, $\alpha_4$ and $\beta$ are arbitrary scalars. Moreover 
this manifold is neither Roter type nor generalized Roter type manifold but 
$S\wedge S$, $S\wedge S^2$, $g\wedge S$, $g\wedge S^2$, $g\wedge g$ and $S^2\wedge S^2$ satisfy the following dependency conditions:
$S\wedge S - \frac{\kappa}{2}\, g\wedge S + \frac{\kappa^2}{4}\, g\wedge g = 0$ and
\beb
& &
L_1\, S\wedge S + L_2\, S\wedge S^2 + L_3\, g\wedge S + L_4\, g\wedge S^2\\
&-& \left(\frac{3}{2} L_3 \kappa + \frac{1}{2} L_1 \kappa^2 + \frac{1}{4} L_4 \kappa^2 + \frac{1}{4} L_2 \kappa^3\right)g\wedge g
+ \left(\frac{8}{\kappa^2} L_3 + \frac{4}{\kappa^2} L_1 - \frac{4}{\kappa^2} L_4\right)S^2\wedge S^2 = 0,
\eeb
where $L_1$, $L_2$, $L_3$ and $L_4$ are arbitrary scalars. Again $(M,g)$ do not satisfy (\ref{b1}), (\ref{b2}) and (\ref{b3}) 
for $T= C, P, K, conh(R)$ except (\ref{b1}) for $T=K$, and also (\ref{b2}), (\ref{b4}) does not hold for $T=R$, $Z = S$ respectively. So the manifold is neither of recurrent curvature $2$-form for $R, C, P, K, conh(R)$ nor of recurrent Ricci $1$-form.\\
\noindent
\textbf{Conclusions.} The study established a bridge between the notions of Deszcz pseudosymmetry with Chaki pseudosymmetry 
as well as Tam$\acute{\mbox{a}}$ssy and Binh's weak symmetry. 
The reduced defining condition of a weakly Ricci symmetric manifold is obtained. 
The independent existence of both the notions of pseudosymmetry be shown by non-trivial examples 
along with the existence of generalized Roter type manifolds. Also the existence 
of the manifolds with recurrent projective, concircular, 
conharmonic curvature $2$-form and recurrent Ricci $1$-form are ensured by concrete examples with various new metrics.\\
\noindent
\textbf{Acknowledgments.} 
The second and fourth named authors are supported by a grant of the Wroc\l aw University of Environmental 
and Life Sciences, Poland. 
The fifth named authors gratefully acknowledges to CSIR, 
New Delhi (File No. 09/025 (0194)/2010-EMR-I) for the financial assistance. 
All the algebraic computations of Section 5 are performed by a program in Wolfram Mathematica,
as well as SymPy and Maxima packages for symbolic calculation.

\indent\\
\noindent \emph{Authors' addresses:}\vspace{-0.36in}
%
\end{document}